\documentclass[a4paper,3p]{elsarticle}
\usepackage[OT4]{fontenc}
\usepackage{amssymb}
\usepackage{amsmath}
\usepackage[all,2cell]{xy}
\UseAllTwocells

\newcommand{\cat}[1]{\ensm{\textup{\bf #1}}}
\newcommand{\catit}[1]{\ensm{\cat{#1}_{it}}}
\newcommand{\algo}{\ensm{\cat{Alg}(\Omega)}}
\newcommand{\latcb}{\ensm{\cat{Lat}_{\bot}}}
\newcommand{\semilatc}[1]{\ensm{\cat{CSLat}(#1)}}
\newcommand{\semilat}{\ensm{\cat{CSLat}(\bigvee)}}
\newcommand{\bornt}[1]{\ensm{\cat{Born}(#1)}}
\newcommand{\lborn}[1]{\ensm{{#1}\text{-}\cat{Born}}}
\newcommand{\lbornsys}[1]{\ensm{{#1}\text{-}\cat{BornSys}}}
\newcommand{\bornsyst}[1]{\ensm{\cat{BornSys}(#1)}}
\newcommand{\set}{\ensm{\cat{Set}}}

\newcommand{\arw}[1]{\ensm{\xrightarrow{#1}}}
\newcommand{\incl}[3]{\ensm{\xymatrix{{#1\,}\ar@{^{(}->}[r]^-{#2} & {#3}}}}

\newcommand{\sqed}{\hfill{\vrule width 3pt height 3pt depth 0pt}}

\newcommand{\ensm}[1]{\ensuremath{#1}}
\newcommand{\mcal}[1]{\ensm{\mathcal{#1}}}
\newcommand{\msf}[1]{\ensm{\mathsf{#1}}}
\newcommand{\seq}{\ensm{\subseteq}}
\newcommand{\prm}[1]{\ensm{{#1}^{\prime}}}
\newcommand{\leqs}{\leqslant}
\newcommand{\opm}[1]{\ensm{{#1}^{op}}}
\newcommand{\gnrt}[1]{\ensm{\langle{#1}\rangle}}

\journal{Fuzzy Sets and Systems}

\begin{document}

\begin{frontmatter}

\title{Categorical foundations of variety-based bornology\tnoteref{GACR}}
\tnotetext[GACR]{The authors gratefully acknowledge the support of Grant Agency of Czech Republic (GA\v{C}R) and
Austrian Science Fund (FWF) within bilateral project No. I 1923-N25 ``New Perspectives on Residuated Posets".}

\author[Paseka]{Jan Paseka}
\address[Paseka]{Department of Mathematics and Statistics, Faculty of Science, Masaryk University\\
                 Kotlarska 2, 611 37 Brno, Czech Republic}
\ead{paseka@math.muni.cz}

\author[Sergey]{Sergey A. Solovyov}
\address[Sergey]{Institute of Mathematics, Faculty of Mechanical Engineering, Brno University of Technology\\
                 Technicka 2896/2, 616 69 Brno, Czech Republic}
\ead{solovjovs@fme.vutbr.cz}

\begin{abstract}
 Following the concept of topological theory of S.~E.~Rodabaugh, this paper introduces a new approach to (lattice-valued) bornology, which is based in bornological theories, and which is called variety-based bornology. In particular, motivated by the notion of topological system of S.~Vickers, we introduce the concept of variety-based bornological system, and show that the category of variety-based bornological spaces is isomorphic to a full reflective subcategory of the category of variety-based bornological systems.
\end{abstract}


\end{frontmatter}


\newtheorem{thm}{Theorem}
\newtheorem{prop}[thm]{Proposition}
\newtheorem{lem}[thm]{Lemma}
\newtheorem{cor}[thm]{Corollary}
\newdefinition{defn}[thm]{Definition}
\newdefinition{exmp}[thm]{Example}
\newdefinition{rem}[thm]{Remark}
\newdefinition{prob}[thm]{Problem}
\newproof{pf}{Proof}


\section{Introduction}


\par
The theory of bornological spaces (which takes its origin in the axiomatisation of the notion of boundedness of
S.-T.~Hu~\cite{Hu1949,Hu1966}) has already found numerous applications in different branches of mathematics. For example, the main ideas of modern functional analysis are those of locally convex topology and convex bornology~\cite{Hogbe-Nlend1977}. Additionally (as a link to physics), one could mention the notion of Hausdorff dimension in convex bornological spaces, motivated by the study of the complexity of strange attractors~\cite{Almeida2002}.

\par
In 2011, M.~Abel and A.~\v{S}ostak~\cite{Abel2011} introduced the concept of lattice-valued (but fixed-basis, in the sense
of~\cite{Hohle1999a}) bornological space, making thereby the first steps towards the theory of lattice-valued bornology. In a series of papers~\cite{Paseka2015,Pasekaa,Paseka2014}, the present authors made further steps in this direction, taking their inspiration in the well-developed theory of lattice-valued topology. In particular, we presented a variable-basis analogue (in the sense of~\cite{Rodabaugh1999a}) of the concept of M.~Abel and A.~\v{S}ostak as well as introduced lattice-valued variable-basis bornological vector spaces (following the pattern of the fuzzy topological vector spaces of~\cite{Katsaras1981,Katsaras1984a}); found the necessary and sufficient condition for the category of lattice-valued bornological spaces to be topological (we notice that all the currently used categories for lattice-valued topology are topological~\cite{Rodabaugh2007}); showed that the category of strict (in the sense of~\cite{Abel2011}) fixed-basis lattice-valued bornological spaces is a quasitopos (we notice again that the categories for lattice-valued topology fail to have this convenient property); and also provided a bornological analogue of the notion of topological system of S.~Vickers~\cite{Vickers1989} (introduced as a common setting for both point-set and point-free topology), starting thus the theory of point-free bornology.

\par
In 2012, S.~Solovyov~\cite{Solovyov2012} presented a new framework for doing lattice-valued topology, namely, variety-based topology. This new setting was induced by an attempt of S.~E.~Rodabaugh~\cite{Rodabaugh2007}, to provide a common framework for the existing categories for lattice-valued topology through the concept of powerset theory (we notice that there also exists the notion of topological theory of  O.~Wyler~\cite{Wyler1971a,Wyler1971}, which is discussed in~\cite{Adamek2009}, and which is compared with the powerset theories of S.~E.~Rodabaugh in~\cite{Rodabaugh2007}). Unlike S.~E.~Rodabaugh, however, who tried to find the basic algebraic structure for lattice-valued topology, arriving thus at the new concept of semi-quantale, S.~Solovyov decided to allow for an arbitrary algebraic structure, thus arriving at varieties of algebras. Briefly speaking, noticing that most of the currently popular approaches to lattice-valued topology are based in powersets of the form $L^X$, where $X$ is a set and $L$ is a complete lattice (possibly, with some additional axioms and/or algebraic operations), S.~Solovyov decided to base his topology in powersets of the form $A^X$, in which the lattice $L$ is replaced with an algebra $A$ from an arbitrary variety (in the extended sense, i.e., allowing for a class of not necessarily finitary operations as in, e.g.,~\cite{Banaschewski1976,Richter1979}). It appeared that such a general variety-based approach is still capable of producing some convenient results, e.g., all the categories of variety-based topological spaces are topological, and, moreover, one has a good concept of topological system, the category of which contains the category of variety-based topological spaces as a full (regular mono)-coreflective subcategory. Even more, variety-based approach to systems incorporates (apart from the lattice-valued extension of topological systems of J.~T.~Denniston, A.~Melton and S.~E.~Rodabaugh~\cite{Denniston2009,Denniston2012}), additionally, state property systems of D.~Aerts~\cite{Aerts1981,Aerts1999a,Aerts2002} (which serve as the basic mathematical structure in the Geneva-Brussels approach to foundations of physics) as well as Chu spaces of~\cite{Pratt1999} (which eventually are nothing else than many-valued formal contexts of Formal Concept Analysis~\cite{Ganter1996}; we notice, however, that Chu spaces were introduced in a more general form by P.-H.~Chu in~\cite{Barr1979}).

\par
The purpose of this paper is to provide a variety-based setting for lattice-valued bornology. In particular, we introduce both variety-based bornological spaces and systems, and show that the category of the former is isomorphic to a full reflective (unlike coreflective in the topological case) subcategory of the latter. Both spaces and systems rely on powerset theories, which, however, are radically different from those of topology. More precisely, one of the main difficulties for introducing a variety-based approach to bornology lies in the fact that while classical topologies are just subframes~\cite{Johnstone1982} of powersets, bornologies are lattice ideals of powersets. Thus, while the switch from subframes to subalgebras is immediate, one needs a good concept of ideal in universal algebras, to switch from lattice ideals to algebra ideals. In this paper, we rely on the concept of algebra ideal of~\cite{Gumm1984} (which eventually goes back to an earlier paper of A.~Ursini). We also notice that motivated by, e.g.,~\cite{Solovyov2013a}, we distinguish between powerset theories and bornological theories. The main reason for this in the topological setting of~\cite{Solovyov2013a} is the fact that while powersets are Boolean algebras, topologies are just subframes, i.e., when dealing with topologies, one ``forgets" a part of the available algebraic structure.

\par
If properly developed, the theory of variety-based bornology could find an application in cancer-related research. More precisely, at the end of Section~7 in~\cite{Hermann2015}, its authors mention that a ``systematic study of continuous spectrum of fractal dimension can put more light on several fractal organisms/objects observed in tissues of cancer patients". Moreover, it appears that  ``in practical applications, we can meet a bornological space instead of a metric one"~\cite{Hermann2015}. Our papers~\cite{Pasekaa,Paseka2014,Paseka2015} together with the present one, could provide a starting point for building a convenient framework for dimension theory for variety-based bornological spaces (e.g., for the already mentioned concept of Hausdorff dimension for bornological spaces of~\cite{Almeida2002}).

\par
This manuscript is based in both category theory and universal algebra, relying more on the former. The necessary categorical background can be found in~\cite{Adamek2009,Herrlich2007}. For algebraic notions, we recommend~\cite{Cohn1981,Gratzer2008,Richter1979}. Although we tried to make the paper as much self-contained as possible, it is expected from the reader to be acquainted with basic concepts of category theory, e.g., with that of reflective subcategory.


\section{Algebraic and categorical preliminaries}


\par
For convenience of the reader, this section briefly recalls those algebraic and categorical concepts, which are necessary for the understanding of this paper.

\par
Variety-based approach to bornology relies on the notion of universal algebra (shortened to algebra), which is thought to be a set, equipped with a family of operations, satisfying certain identities. The theory of universal algebra calls a class of finitary algebras (induced by a set of finitary operations), closed under the formation of homomorphic images, subalgebras and direct products, a variety. In view of the structures common in, e.g., lattice-valued topology (where set-theoretic unions are replaced with joins), we consider infinitary algebraic theories, extending the approach of varieties to cover our needs (cf. \cite{Manes1976,Richter1979}).

\begin{defn}
 \label{defn:1}
 Let $\Omega=(n_\lambda)_{\lambda\in\Lambda}$ be a (possibly
 proper or empty) class of cardinal numbers. An \emph{$\Omega$-algebra} is a pair $(A,(\omega^A_\lambda)_{\lambda\in\Lambda})$, which comprises a set $A$ and a family of maps $A^{n_\lambda}\arw{\omega^A_\lambda}A$ (\emph{$n_\lambda$-ary primitive operations} on $A$). An \emph{$\Omega$-homomorphism} $(A,(\omega^A_\lambda)_{\lambda\in
 \Lambda})\arw{\varphi}(B,(\omega^B_\lambda)_{\lambda\in\Lambda})$ is a map $A\arw{\varphi}B$, which makes the diagram
 $$
   \xymatrix@1{A^{n_\lambda} \ar[d]_-{\omega^A_\lambda}
               \ar[rr]^-{\varphi^{n_\lambda}} && B^{n_\lambda}
               \ar[d]^-{\omega^B_\lambda} \\
               A \ar[rr]_-{\varphi} && B\\}
 $$
 commute for every $\lambda\in\Lambda$. \algo\ is the construct of $\Omega$-algebras and $\Omega$-homomorphisms.
 \sqed
\end{defn}

\par
Every concrete category of this paper has the underlying functor $|-|$ to the respective ground category (e.g., the category \cat{Set} of sets and maps in case of constructs), the latter mentioned explicitly in each case.

\begin{defn}
 \label{defn:2}
 Let \mcal{M} (resp. \mcal{E}) be the class of $\Omega$-homomorphisms with injective (resp. surjective) underlying maps. A \emph{variety of $\Omega$-algebras} is a full subcategory of \algo, which is closed under the formation of products, \mcal{M}-subobjects (subalgebras) and \mcal{E}-quotients (homomorphic images). The objects (resp. morphisms) of a variety are called \emph{algebras} (resp. \emph{homomorphisms}).
 \sqed
\end{defn}

\begin{defn}
 \label{defn:3}
 Given a variety \cat{A}, a \emph{reduct} of \cat{A} is a pair $(\|-\|,\cat{B})$, where \cat{B} is a variety such that $\Omega_{\cat{B}}\seq\Omega_{\cat{A}}$ (every primitive operation of \cat{B} is a primitive operation of \cat{A}), and $\cat{A}\arw{\|-\|}\cat{B}$ is a concrete functor. The pair $(\cat{A},\|-\|)$ is called then an \emph{extension} of \cat{B}.
 \sqed
\end{defn}

\par
The following example illustrates the concept of variety with several well-known constructs.

\begin{exmp}
 \label{exmp:1}
 \hfill\par
 \begin{enumerate}[(1)]
  \item \latcb\ is the variety of \emph{lattices with the smallest element}, i.e., partially ordered sets (posets), which have binary
        $\wedge$ and finite (including the empty) $\vee$.
  \item \semilatc{\Xi} is the variety of \emph{$\Xi$-semilattices}, i.e., posets, which have arbitrary $\Xi\in\{\bigwedge,\bigvee\}$.
  \item \cat{SQuant} is the variety of \emph{semi-quantales}, i.e., $\bigvee$-semilattices, equipped with a binary
        operation $\otimes$~\cite{Rodabaugh2007}.
  \item \cat{SFrm} is the variety of \emph{semi-frames}, i.e., $\bigvee$-semilattices, with singled out finite meets~\cite{Rodabaugh2007}.
  \item \cat{QSFrm} is the variety of \emph{quantic semi-frames}, i.e., semi-frames, which are also semi-quantales, satisfying the
        conditions $a\otimes\top_A=\top_A\otimes a=a$ (\emph{strict two-sidedness}) and $a\otimes\bot_A=\bot_A\otimes a=\bot_A$ for every $a\in A$.
  \item \cat{Frm} is the variety of \emph{frames}, i.e., semi-frames $A$, which additionally satisfy the
        distributivity condition $a\wedge(\bigvee S)=\bigvee_{s\in S}(a\wedge s)$ for every $a\in A$ and every $S\seq A$~\cite{Johnstone1982}.
  \item \cat{CLat} is the variety of \emph{complete lattices}, i.e., posets $A$, which are both $\bigvee$- and $\bigwedge$-semilattices.
  \item \cat{CBAlg} is the variety of \emph{complete Boolean algebras}, i.e., complete lattices $A$ such that $a\wedge(b\vee
        c)=(a\wedge b)\vee(a\wedge c)$ for every $a$, $b$, $c\in A$, equipped with a unary operation $A\arw{(-)^{\ast}}A$ such that $a\vee a^{\ast}=\top_A$ and $a\wedge a^{\ast}=\bot_A$ for every $a\in A$, where $\top_A$ (resp. $\bot_A$) is the largest (resp. smallest) element of $A$.
        \sqed
 \end{enumerate}
\end{exmp}

\par
We notice that, e.g., \cat{CLat} is a reduct of \cat{CBAlg}, and both \semilatc{\Xi} and \latcb\ are reducts of \cat{CLat}. In the next step, following~\cite{Gumm1984}, we recall the concept of ideal in universal algebras. We underline, however, immediately that while~\cite{Gumm1984} works with classical universal algebras (which have a set of finitary operations), we adapt the respective definition of ideal to our case of generalized varieties. Moreover, for the sake of convenience (and following~\cite{Gumm1984}), given a set $X$ and a cardinal number $n$ (not necessarily finite), the elements of the powerset $X^n$ will be denoted $\vec{x}$, or sometimes (for better clarity) $\vec{x}_{i\in n}$. Additionally, given an element $x_0\in X$, $\vec{x_0}$ will denote the element of $X^n$, every component of which is $x_0$.

\begin{defn}
 \label{defn:4}
 Let \cat{A} be a variety of algebras, the elements of which have a 0-ary operation $\bot$.
 \begin{enumerate}[(1)]
  \item A term $p(\vec{x},\vec{y})$ is an \emph{(\cat{A}-)ideal term in $\vec{y}$} if $p(\vec{x},\vec{\bot})=\bot$ is an identity in
        \cat{A}.
  \item A non-empty subset \mcal{I} of an \cat{A}-algebra $A$ is called an \emph{(\cat{A}-)ideal} provided that for every ideal term
        $p(\vec{x},\vec{y})$ in $\vec{y}$, $\vec{a}\in A^n$, $\vec{i}\in\mcal{I}^m$, it follows that $p(\vec{a},\vec{i})\in\mcal{I}$.
        \sqed
 \end{enumerate}
\end{defn}

\par
Every ideal contains the constant $\bot$, since the (constant) term $\bot$ is an ideal term. Moreover, in the varieties \cat{Rng} (rings) and
\latcb, ideals in the sense of Definition~\ref{defn:4} coincide with their respective classical analogues. To use it later on in the paper, we recall the definition of lattice ideals.

\begin{defn}
 \label{defn:5}
 Given a \latcb-object $L$ with the smallest element $\bot_L$, a subset $\mcal{I}\seq L$ is said to be an \emph{ideal} of $L$ provided that
 the following conditions are fulfilled:
 \begin{enumerate}[(1)]
  \item for every finite subset $S\seq L$, $\bigvee S\in\mcal{I}$;
  \item for every $a\in L$ and every $i\in I$, $a\wedge i\in\mcal{I}$.
        \sqed
 \end{enumerate}
\end{defn}

\par
Definition~\ref{defn:5}\,(1) implies that every lattice ideal is non-empty, i.e., $\bot_L\in\mcal{I}$. Moreover, Definition~\ref{defn:5} suggests two possible (but not the only ones) ideal terms, i.e., $\bigvee\vec{y}$ (the term takes finite sequences and does not depend on $\vec{x}$) and $x\wedge y$. All the other lattice ideal terms are generated by the just mentioned two. For example, one can build ideal terms of the form $(x_1\wedge y_1)\vee(x_2\wedge y_2)$ or $(x\wedge y_1)\vee y_2$.

\par
To use it later on in the paper, we recall a convenient property of ideals in universal algebras from~\cite{Gumm1984}.

\begin{prop}[\cite{Gumm1984}]
 \label{prop:1}
 Let \cat{A} be a variety with a 0-ary operation, and let $A$ be an \cat{A}-algebra.
 \begin{enumerate}[(1)]
  \item The intersection of a family of ideals of $A$ is an ideal.
  \item Given a subset $S\seq A$, the ideal generated by $S$ (denoted $\gnrt{S}_{\cat{A}}$) consists of all $p(\vec{a},\vec{s})$, where
        $p(\vec{x},\vec{y})$ is an ideal term in $\vec{y}$, $\vec{a}\in A^n$, and $\vec{s}\in S^m$.
 \end{enumerate}
\end{prop}

\par
In the last step, we introduce a particular construction, which is related to ideals in universal algebras.

\begin{defn}
 \label{defn:6}
 Given a variety of algebras \cat{A} with a 0-ary operation, \catit{A} is the construct, whose objects are \cat{A}-ideals, and whose morphisms are maps $\mcal{I}_1\arw{\varphi}\mcal{I}_2$, which preserve the ideal terms as follows:
 \begin{enumerate}[(1)]
  \item for every ideal term $p(\vec{y})$ in $\vec{y}$, $\varphi(p(\vec{i}))=p(\varphi(\vec{i}))$ for every $\vec{i}\in\mcal{I}_1^n$;
  \item for every ideal term $p(\vec{x},\vec{y})$ in $\vec{y}$, there exist an ideal term
        $t(\vec{x}\cup\vec{z},\vec{y})$ in $\vec{y}$ (in which the notation ``$\vec{x}\cup\vec{z}$" means that one adds more variables to $\vec{x}$, e.g., taking $x_1,x_2,z_1,z_2$ instead of $x_1,x_2$, but not to $\vec{y}$), and a family of ideal terms $\overrightarrow{s_u(\vec{v_u},\vec{w_u})}_{u\in\vec{z}}$ in $\vec{w_u}$ such that $\varphi(p(\vec{a},\vec{i}))=t(\varphi(\vec{a})\cup \varphi(\overrightarrow{s_u(\vec{\prm{a}_u},\vec{\prm{i}_u})}_{u\in\vec{z}}),\varphi(\vec{i}))$ for every $\vec{a}\in\mcal{I}_1^n$, $\vec{i}\in\mcal{I}_1^m$, where $\vec{\prm{a}_u}\seq\vec{a}$, $\vec{\prm{i}_u}\seq\vec{i}$.
        \sqed
 \end{enumerate}
\end{defn}

\par
The main result, which could be easily obtained with the help of Definition~\ref{defn:6} and which will be used throughout the paper, can be stated now as follows.

\begin{prop}
 \label{prop:1.1}
 Given an \catit{A}-morphism $\mcal{I}_1\arw{\varphi}\mcal{I}_2$, for every ideal $\mcal{J}_2\seq\mcal{I}_2$, the set $\mcal{J}_1=\{a\in\mcal{I}_1\,|\,\varphi(a)\in\mcal{J}_2\}$ is an ideal.
\end{prop}
\begin{pf}
 It follows directly from the two items of Definition~\ref{defn:6} that $\mcal{J}_1$ is closed under the ideal terms.
 \qed
\end{pf}

\par
For convenience of the reader, we notice that, eventually, we could have replaced the characterization of \catit{A}-morphisms in Definition~\ref{defn:6} with the one in Proposition~\ref{prop:1.1}, namely, with the property that the preimage of an ideal is an ideal. The main intuition for Definition~\ref{defn:6}, however, will be provided while discussing the topic of powerset theories in the next section. At this place, we just mention the following easy result.

\begin{prop}
 \label{prop:2}
 Given a variety \cat{A} with a 0-ary operation, there exists a non-full embedding \incl{\cat{A}}{E}{\catit{A},} which is defined by $E(A_1\arw{\varphi}A_2)=A_1\arw{\varphi}A_2$.
\end{prop}
\begin{pf}
 Given an \cat{A}-homomorphism $A_1\arw{\varphi}A_2$ and an ideal term $p(\vec{x},\vec{y})$ in $\vec{y}$, it follows that $\varphi(p(\vec{a},\vec{i}))=p(\varphi(\vec{a}),\varphi(\vec{i}))$ for every $\vec{a}\in A_1^n$, $\vec{i}\in A_1^m$.
 \qed
\end{pf}

\par
For the sake of convenience, later on, we will not distinguish between \cat{A} and its image in \catit{A} under $E$.

\par
To conclude the preliminaries, we recall from~\cite{Adamek2009} the definition of reflective subcategory (the dual case of coreflective subcategories is left to the reader).

\begin{defn}
 \label{den:6.1}
 Let \cat{A} be a subcategory of \cat{B}, and let $B$ be a \cat{B}-object.
 \begin{enumerate}[(1)]
  \item An \emph{\cat{A}-reflection} (or \emph{\cat{A}-reflection arrow}) for $B$ is a \cat{B}-morphism $B\arw{r}A$ from $B$ to an
        \cat{A}-object $A$ with the following universal property: for every \cat{B}-morphism $B\arw{f}\prm{A}$ from $B$ into some \cat{A}-object \prm{A}, there exists a unique \cat{A}-morphism $A\arw{\prm{f}}\prm{A}$ such that the triangle
        $$
          \xymatrix{B \ar[r]^-{r} \ar[rd]_-{f} & A \ar@{.>}[d]^-{\prm{f}}\\
                    & \prm{A}}
        $$
        commutes.
  \item \cat{A} is called a \emph{reflective subcategory} of \cat{B} provided that each \cat{B}-object has an \cat{A}-reflection.
        \sqed
 \end{enumerate}
\end{defn}

\par
To give the reader more intuition, we provide one example of reflective subcategories.

\begin{exmp}
 \label{exmp:1.1}
 The category \cat{Mon} of monoids is a reflective subcategory of the category \cat{Sgr} of semigroups. Given a semigroup $(X,\circ)$, the extension \incl{(X,\circ)}{}{(X\cup\{e\},\hat{\circ},e),} which is obtained by adding a unit element $e\not\in X$ of the operation $\hat{\circ}$, is a \cat{Mon}-reflection for $(X,\circ)$.
 \sqed
\end{exmp}


\section{Variety-based bornology}


\par
This section introduces our new approach to bornology, which is based in varieties of algebras. To provide more intuition for the next developments, we begin by recalling the concept of bornological space from~\cite{Hogbe-Nlend1977}.

\par
We notice first that there exists the so-called covariant powerset functor $\set\arw{\mcal{P}}\set$, which is defined by $\mcal{P}(X\arw{f}Y)=\mcal{P}X\arw{\mcal{P}f}\mcal{P}Y$, where $\mcal{P}f(S)=\{f(s)\,|\,s\in S\}$.

\begin{defn}
 \label{defn:7}
 A \emph{bornological space} is a pair $(X,\mcal{B})$, where $X$ is a set and $\mcal{B}$ (a \emph{bornology} on $X$) is a subfamily of $\mcal{P}X$ (the elements of which are called \emph{bounded sets}), which satisfies the following axioms:
 \begin{enumerate}[(1)]
  \item $X=\bigcup\mcal{B}$ $(=\bigcup_{B\in\mcal{B}}B)$;
  \item if $B\in\mcal{B}$ and $D\seq B$, then $D\in\mcal{B}$;
  \item if $\mcal{S}\seq\mcal{B}$ is finite, then $\bigcup\mcal{S}\in\mcal{B}$.
 \end{enumerate}
 Given bornological spaces $(X_1,\mcal{B}_1)$ and $(X_2,\mcal{B}_2)$, a map $X_1\arw{f}X_2$ is called \emph{bounded} provided that $\mcal{P}f(B_1)\in\mcal{B}_2$ for every $B_1\in\mcal{B}_1$. \cat{Born} is the construct of bornological spaces and bounded maps.
 \sqed
\end{defn}

\par
To start a short discussion, we notice that given a map $X\arw{f}Y$, while $\mcal{P}X$, $\mcal{P}Y$ are complete Boolean algebras, the map $\mcal{P}f$ is just a \semilat-homomorphism, i.e., $\mcal{P}(\bigcup\mcal{S})=\bigcup_{S\in\mcal{S}}\mcal{P}(S)$ for every $\mcal{S}\seq\mcal{P}X$. We underline, moreover, that given $S,T\in\mcal{P}X$, $\mcal{P}f(S\bigcap T)=
(\mcal{P}f(S)\bigcap\mcal{P}f(S\bigcap T))\bigcap\mcal{P}f(T)$, which takes us back to Definition~\ref{defn:6} (recall the discussion after Definition~\ref{defn:5} on ideal terms for lattices, and notice that $(u\wedge v)\wedge w$ is a lattice ideal term in $w$). As a consequence, we could have considered (verification is easy) the functor $\set\arw{\mcal{P}}\catit{CBAlg}$ instead of just $\set\arw{\mcal{P}}\set$.

\par
To continue, we notice that items~(2), (3) of Definition~\ref{defn:7} ensure that every bornology on a set $X$ is a lattice ideal of $\mcal{P}X$, but not a complete lattice ideal of $\mcal{P}X$ in the sense of Definition~\ref{defn:4} (no closure under arbitrary joins). Thus, bornologies on sets are the elements of the variety \catit{\latcb} and not \catit{CBAlg}.

\par
Lastly, we notice that the first item of Definition~\ref{defn:7} just says that the \catit{CBAlg}-ideal of $\mcal{P}X$, generated by $\mcal{B}$, contains the largest element of $\mcal{P}X$ (which implies then that $\gnrt{\mcal{B}}_{\catit{CBAlg}}$ is the whole $\mcal{P}X$).

\par
The above discussion gives rise to the next variety-based analogue of the classical approach to bornology.


\subsection{Variety-based powerset and bornological theories}


\par
From now on, we assume that every variety \cat{A} has 0-ary operations $\bot$ and $\top$, and every our employed reduct \cat{B} of \cat{A} preserves the operation $\bot$ (but not necessarily $\top$). One should also recall Proposition~\ref{prop:2}.

\begin{defn}
 \label{defn:8}
 A \emph{powerset theory} in a category \cat{X} (\emph{ground category} of the theory) is a functor $\cat{X}\arw{P}\catit{A}$ to the category of ideals of a variety \cat{A}, with the additional property that $PX$ lies in \cat{A} for every $X\in\cat{X}$.
 \sqed
\end{defn}

\par
The following example illustrates the just introduced notion.

\begin{exmp}
 \label{exmp:2}
 \hfill\par
 \begin{enumerate}[(1)]
  \item The classical covariant powerset functor $\set\arw{\mcal{P}}\set$ provides a powerset theory $\set\arw{\mcal{P}}\catit{CBAlg}$.
  \item Given a $\bigvee$-semilattice $L$, there is a functor $\set\arw{\mcal{P}_L}\semilat$, $\mcal{P}_L(X\arw{f}Y)=
        L^X\arw{\mcal{P}_Lf}L^Y$, where $(\mcal{P}_Lf(\alpha))(y)=\bigvee_{f(x)=y}\alpha(x)$~\cite{Rodabaugh2007}. For every ``suitable" extension \cat{A} of \semilat\ (i.e., \cat{SQuant}, \cat{SFrm}, \cat{QSFrm}, \cat{Frm}, \cat{CLat}, \cat{CBAlg}), one gets then a powerset theory $\set\arw{\mcal{P}_A}\catit{A}$ (which will be denoted $\mcal{P}_A^{\cat{A}}$). In particular, the two-element Boolean algebra $\msf{2}=\{\bot,\top\}$ provides the powerset theory from the previous item. For convenience of the reader, we also emphasize here that in case of the variety \cat{QSFrm}, both $x\otimes y$ and $y\otimes x$ are ideal terms in $y$, for which one can easily show that $\mcal{P}_Lf(\alpha\otimes\beta)=(\mcal{P}_Lf(\alpha)\wedge\mcal{P}_Lf(\alpha\otimes\beta))\wedge\mcal{P}_L(\beta)$ and use the fact that $x\wedge y$ is an ideal term in $y$ (notice that, by the definition of \cat{QSFrm}, $(\mcal{P}_Lf(\alpha\otimes\beta))(y)=\bigvee_{f(x)=y}(\alpha(x)\otimes\beta(x))\leqs\bigvee_{f(x)=y}
        (\alpha(x)\otimes\top_A)=\bigvee_{f(x)=y}\alpha(x)=(\mcal{P}_Lf(\alpha))(y)$). In case of the variety \cat{SQuant}, however, $x\otimes y$ (also $y\otimes x$) is not an ideal term in $y$.
  \item Given a subcategory \cat{S} of \semilat, there is a functor $\set\times\cat{S}\arw{\mcal{P}_{\cat{S}}}\semilat$, which is defined
        by $\mcal{P}_{\cat{S}}((X,L)\arw{(f,\varphi)}(Y,M))=L^X\arw{\mcal{P}_{\cat{S}}(f,\varphi)}M^Y$, where $(\mcal{P}_{\cat{S}}f(\alpha))(y)=\bigvee_{f(x)=y}\varphi\circ\alpha(x)$~\cite{Rodabaugh2007}. For every ``suitable" extension \cat{A} of \semilat\ (i.e., \cat{SQuant}, \cat{SFrm}, \cat{QSFrm}, \cat{Frm}, \cat{CLat}, \cat{CBAlg}), one obtains then a powerset theory $\set\times\cat{S}\arw{\mcal{P}_{\cat{S}}}\catit{A}$ (which will be denoted $\mcal{P}_{\cat{S}}^{\cat{A}}$). If \cat{S} is the subcategory of the form $L\arw{1_L}L$, then one gets the powerset theory of the previous item.~\sqed
 \end{enumerate}
\end{exmp}

\par
For convenience of the reader, we notice the crucial difference between bornological and topological powerset
theories~\cite{Rodabaugh2007,Solovyov2013a}: while the former extend the covariant powerset functor, the latter extend the so-called contravariant powerset functor of the form $\set\arw{\mcal{Q}}\opm{\set}$, which is given by $\mcal{Q}(X\arw{f}Y)=\mcal{P}X\arw{\mcal{Q}f}\mcal{P}Y$,
$\opm{(\mcal{Q}f)}(S)=\{x\in X\,|\,f(x)\in S\}$. Since \mcal{Q} has, eventually, the form $\set\arw{Q}\opm{\cat{CBAlg}}$, the general form of a (topological) powerset theory is a functor $\cat{X}\arw{P}\opm{\cat{A}}$, where \cat{A} is a variety of algebras.

\par
We are ready to define bornological theories, which "forget" certain algebraic structure of "powersets".

\begin{defn}
 \label{defn:9}
 Given a powerset theory $\cat{X}\arw{P}\catit{A}$ and a reduct \cat{B} of \cat{A}, a \emph{bornological theory} in \cat{X} induced by
 the pair $(P,\cat{B})$ is the functor $\cat{X}\arw{T}\catit{B}$, which is defined as the composition $\cat{X}\arw{P}\catit{A}\arw{\|-\|_{it}}\catit{B}$.
 \sqed
\end{defn}

\par
As an example of bornological theories, we mention the classical one $\set\arw{T}\catit{\latcb}$, which is induced by the powerset theory $\set\arw{\mcal{P}}\catit{CBAlg}$ and the reduct \latcb\ of \cat{CBAlg}.


\subsection{Variety-based bornological spaces}


\par
With bornological theories in hand, this subsection introduces variety-based bornological spaces. The reader should recall the existence of the 0-ary operation $\top$ in every variety \cat{A} we use.

\begin{defn}
 \label{defn:10}
  Let $T$ be a bornological theory in a category \cat{X}. \bornt{T} is the concrete category over \cat{X}, whose objects
  (\emph{$T$-bornological spaces} or \emph{$T$-spaces}) are pairs $(X,\tau)$, where $X$ is an \cat{X}-object and $\tau$ (\emph{$T$-bornology} on $X$) is an ideal of $TX$ such that $\top\in\gnrt{\tau}_{\cat{A}}$, and whose morphisms (\emph{$T$-bounded \cat{X}-morphisms}) $(X_1,\tau_1)\arw{f}(X_2,\tau_2)$ are \cat{X}-morphisms $X_1\arw{f}X_2$ such that $Tf(\alpha)\in\tau_2$ for every $\alpha\in\tau_1$ (\emph{$T$-boundedness}).
 \sqed
\end{defn}

\par
The next example illustrates the just introduced concept of variety-based bornology.

\begin{exmp}
 \label{exmp:3}
 \hfill\par
 \begin{enumerate}[(1)]
  \item \bornt{(\mcal{P},\latcb)} is isomorphic to the category \cat{Born} of bornological spaces and bounded maps
        of~\cite{Hogbe-Nlend1977}.
  \item \bornt{(\mcal{P}_{L}^{\cat{CLat}},\latcb)} is isomorphic to the category \lborn{L} of~\cite{Abel2011}.
  \item \bornt{(\mcal{P}_{\cat{S}}^{\cat{CLat}},\latcb)} is isomorphic to the category \lborn{\cat{L}} of~\cite{Paseka2015}.
        \sqed
 \end{enumerate}
\end{exmp}

\par
We have already mentioned that every category of variety-based topological spaces is topological over its ground category~\cite{Solovyov2012,Solovyov2013a}. In~\cite{Paseka2015,Pasekaa}, we showed that not every category of lattice-valued bornological spaces is topological, providing the necessary and sufficient conditions (on the respective lattices) for this property. Below, we show sufficient conditions for the category \bornt{T} to be topological over \cat{X}.

\par
To provide more intuition for the subsequent general requirements, we start with a motivating discussion, which recalls some well-known (and simple) properties of complete lattices and their induced powersets.

\par
First, given a complete lattice $L$ and a lattice ideal $\mcal{I}\seq L$ such that $\top_L\in\gnrt{\mcal{I}}_{\cat{CLat}}$, there exists a subset $S\seq\mcal{I}$ with $\bigvee S=\top_L$ (cf. Proposition~\ref{prop:1}). We can consider the operation $\bigvee$ as a term $t(\vec{x})$ of no fixed arity.

\par
Second, given complete lattices $L$, $M$, every $\bigvee$-semilattice morphism $L\arw{\varphi}M$ has a right adjoint (in the sense of posets) map $|L|\arw{\varphi^{\vdash}}|M|$, which is defined by $\varphi^{\vdash}(b)=\bigvee\{a\in L\,|\,\varphi(a)\leqs b\}$, and which, additionally, is $\bigwedge$-preserving. In particular, $\varphi\circ\varphi^{\vdash}\leqs 1_{M}$, which implies that given a lattice ideal $\mcal{I}\seq M$, $\varphi\circ\varphi^{\vdash}(i)\in\mcal{I}$ for every $i\in\mcal{I}$. Further, given a map $X\arw{f}Y$, Example~\ref{exmp:2}\,(3) provides a $\bigvee$-preserving map $L^X\arw{\mcal{P}_{\cat{S}}(f,\varphi)}M^Y$, which (by the above formula) has a right adjoint $|M^Y|\arw{(\mcal{P}_{\cat{S}}(f,\varphi))^{\vdash}}|L^X|$ given by $(\mcal{P}_{\cat{S}}(f,\varphi))^{\vdash}(\beta)=\varphi^{\vdash}\circ\beta\circ f$, and which, moreover, is a complete lattice homomorphism provided that $\varphi^{\vdash}$ is $\bigvee$-preserving.

\par
Third, with the notations of the above item, for every subset $\mcal{L}\seq L^X$, there exists the meet $\bigwedge\mcal{L}$ such that $(\mcal{P}_{\cat{S}}(f,\varphi))(\bigwedge\mcal{L})\leqs(\mcal{P}_{\cat{S}}(f,\varphi))(\alpha)$ for every $\alpha\in\mcal{L}$, i.e., $(\mcal{P}_{\cat{S}}(f,\varphi))(\bigwedge\mcal{L})$ belongs to the lattice ideal of $M^Y$ generated by $(\mcal{P}_{\cat{S}}(f,\varphi))(\alpha)$ for every $\alpha\in\mcal{L}$. Additionally, $\bigwedge\mcal{L}=\underline{\top_L}$ provided that $\alpha=\underline{\top_L}$ for every $\alpha\in\mcal{L}$. We can consider the operation $\bigwedge$ as a term $t(\vec{x})$ with no fixed arity.

\par
Fourth, in~\cite{Paseka2015,Pasekaa}, we showed that the category \lborn{L} (recall Example~\ref{exmp:3}\,(2)) is a topological construct iff the complete lattice $L$ is ideally completely distributive at $\top_L$, which is defined as follows.

\begin{defn}
 \label{defn:11}
 A complete lattice $L$ is called \emph{ideally completely distributive at $\top_L$} provided that for every non-empty family $\{S_i\,|\,i\in I\}$ of lattice ideals of $L$, $\bigwedge_{i\in I}(\bigvee S_i)=\top_L$ implies $\bigvee_{h\in H}(\bigwedge_{i\in I}h(i))=\top_L$, where $H$ is the set of all \emph{choice maps} on $\bigcup_{i\in I}S_i$, i.e., maps $I\arw{h}\bigcup_{i\in I}S_i$ such that $h(i)\in S_i$ for every $i\in I$.
 \sqed
\end{defn}

\par
For convenience of the reader, we notice that ideal complete distributivity at $\top$ does not even imply distributivity; conversely, every completely distributive lattice~\cite{Gierz2003} is ideally completely distributive at $\top$, but there exists an infinitely
distributive lattice $L$ (namely, $a\wedge(\bigvee S)=\bigvee_{s\in S}(a\wedge s)$ and $a\vee(\bigwedge S)=\bigwedge_{s\in S}(a\vee s)$ for every $a\in L$ and every $S\seq L$), which is not ideally completely distributive at $\top$~\cite{Paseka2015}.

\par
We introduce now four requirements on our variety-based bornological setting, which correspond to the above-mentioned four items.

\begin{description}
 \item[Req.~1] Variety \cat{A} has an ideal term $p(\vec{x})$ of no fixed arity (for every \cat{A}-algebra $A$ and every $\vec{a}\in
               A^n$, $p(\vec{a})$ is defined) such that given an \cat{A}-algebra $A$ and an ideal $\mcal{I}\seq A$ containing $\top$, $p(\vec{i})=\top$ for some $\vec{i}\in\mcal{I}^n$.
\end{description}

\begin{description}
 \item[Req.~2] Every \catit{A}-morphism $A_1\arw{\varphi}A_2$ has a map $|A_2|\arw{\varphi^{\vdash}}|A_1|$ such that for every
               \cat{B}-ideal $\mcal{I}\seq A_2$ and every $i\in\mcal{I}$, $\varphi\circ\varphi^{\vdash}(i)\in\mcal{I}$. If $PX_1\arw{\varphi}PX_2$ for some \cat{X}-objects $X_1$, $X_2$, then $\varphi^{\vdash}$ has two additional properties: first, $\varphi^{\vdash}(p(\vec{x}))=p(\varphi^{\vdash}(\vec{x}))$ for every $\vec{x}\in(PX_1)^n$ such that $p(\vec{x})=\top$, where $p(\vec{y})$ is the ideal term from Req.~1; and, second, $\varphi^{\vdash}$ preserves $\top$.
\end{description}

\begin{description}
 \item[Req.~3] Variety \cat{A} has a term $t(\vec{x})$ of no fixed arity such that, first, for every \catit{A}-morphism
               $A_1\arw{\varphi}A_2$, $\varphi(t(\vec{a}))\in\gnrt{\varphi(b)}_{\cat{B}}$ for every $b\in\vec{a}$ (in which the notation ``$b\in\vec{a}$" means that $b$ is one of the components of the sequence $\vec{a}$), and, second, $t(\vec{a})=\top$ provided that $b=\top$ for every $b\in\vec{a}$.
\end{description}

\begin{description}
 \item[Req.~4] For every \cat{X}-object $X$, $PX$ is \emph{ideally completely distributive at $\top$}, i.e., given a family $\vec{y}_j\in
               (PX)^{n_j}$ with $j\in J$, $t(\overrightarrow{p(\vec{y}_j)}_{j\in J})=\top$ implies $p(\overrightarrow{t(\overrightarrow{h(j)}_{j\in J})}_{h\in H})=\top$, with $H$ the set of choice maps on $\bigcup_{j\in J}\vec{y_j}$.
\end{description}

\par
The necessary requirements in hand, we can now prove the main results of this subsection.

\begin{thm}
 \label{thm:1}
 If Req.~1~--~4 hold, then the category \bornt{T} is topological over \cat{X}.
\end{thm}
\begin{pf}
 Since the category \bornt{T} is clearly amnestic, by~\cite[Proposition~21.5]{Adamek2009}, it will be enough to show that every structured
 $|-|$-source $(X\arw{f_j}|(X_j,\tau_j)|)_{j\in J}$ has an $|-|$-initial lift. Define $\tau=\{\alpha\in TX\,|\,Tf_j(\alpha)\in\tau_j$ for every $j\in J\}$. We have to show that, first, $\tau$ is an ideal of $TX$, and, second, $\top\in\gnrt{\tau}_{\cat{A}}$.

 \par
 The first claim follows immediately from Propositions~\ref{prop:1}, \ref{prop:1.1}. To show that $\top\in\gnrt{\tau}_{\cat{A}}$, we notice first that for every $j\in J$, by Req.~1, there exists $\vec{y_j}\in\tau_j^{n_j}$ such that $p(\vec{y_j})=\top_j$. Let $H$ be the set of choice functions on $\bigcup_{j\in J}\vec{y_j}$. For every $h\in H$, by Req.~2, 3, we define $\alpha_h=t(\overrightarrow{(Tf_j)^{\vdash}h(j)}_{j\in J})$. To show that $\alpha_h\in\tau$, we notice that given $j_0\in J$, $Tf_{j_0}(\alpha_h)=Tf_{j_0}(t(\overrightarrow{(Tf_j)^{\vdash}h(j)}_{j\in J}))\overset{Req.~3}{\in}\gnrt{Tf_{j_0}\circ(Tf_{j_0})^{\vdash}(h(j_0))}_{\cat{B}}\overset{Req.~2}{\seq}\tau_{j_0}$. As a consequence, we get (by Req.~4) that $t(\overrightarrow{p((Tf_j)^{\vdash}(\vec{y_j}))}_{j\in J})\overset{Req.~2}{=}
 t(\overrightarrow{(Tf_j)^{\vdash}(p(\vec{y_j}))}_{j\in J})=t(\overrightarrow{(Tf_j)^{\vdash}(\top_j)}_{j\in J})\overset{Req.~2}{=}
 t(\vec{\top})\overset{Req.~3}{=}\top$ implies $\top=p(\overrightarrow{t(\overrightarrow{(Tf_j)^{\vdash}(h(j))}_{j\in J})}_{h\in H})=p(\overrightarrow{\alpha_h}_{h\in H})$, i.e., $\top\in\gnrt{\tau}_{\cat{A}}$.
 \qed
\end{pf}

\par
We notice that Theorem~\ref{thm:1} implies, in particular, one part (the sufficiency) of the result of~\cite{Paseka2015} on the topological nature of the category \lborn{\cat{L}} (recall Example~\ref{exmp:3}\,(3)).

\begin{defn}
 \label{defn:12}
 $\cat{L}^{\vdash}$ is the subcategory of \cat{CLat}, the objects of which are complete lattices $L$, which are
 ideally completely distributive at $\top_L$, and whose morphisms $L_1\arw{\psi}L_2$ are such that the map $L_2\arw{\psi^{\vdash}}L_1$ has the following property: $\psi^{\vdash}(\bigvee S)=\bigvee_{\beta\in S}\psi^{\vdash}(\beta)$ for every $S\seq L_2$ with $\bigvee S=\top_{L_2}$.
 \sqed
\end{defn}

\par
From Theorem~\ref{thm:1}, one then immediately gets the next result (notice that validity of Req.~1~--~4 in this particular case served as our main motivation for their introduction).

\begin{thm}
 \label{thm:2}
 If \cat{L} is a subcategory of $\cat{L}^{\vdash}$, then the category \lborn{\cat{L}} is topological.
\end{thm}


\subsection{Variety-based bornological systems}


\par
In this subsection, we introduce a variety-based bornological analogue of the concept of topological system of S.~Vickers~\cite{Vickers1989}. We notice, however, immediately that~\cite{Paseka2015} has already introduced lattice-valued bornological systems, motivated by the notion of lattice-valued topological system of~\cite{Denniston2009,Denniston2012}. It is the purpose of this subsection, to provide a bornological analogue of variety-based topological systems of~\cite{Solovyov2012,Solovyov2013a}.

\begin{defn}
 \label{defn:13}
 Given a bornological theory $\cat{X}\arw{T}\catit{B}$, \bornsyst{T} is the comma category ${(1_{\catit{B}}\downarrow T)}$, concrete over the product category $\cat{X}\times\catit{B}$, whose objects (\emph{$T$-bornological systems} or \emph{$T$-systems}) are triples $(X,\kappa,B)$, which are made by \catit{B}-morphisms $B\arw{\kappa}TX$ such that $\top\in\gnrt{\mcal{P}\kappa(B)}_{\cat{A}}$, and whose morphisms (\emph{$T$-bounded morphisms}) $(X_1,\kappa_1,B_1)\arw{(f,\varphi)}(X_2,\kappa_2,B_2)$ are $\cat{X}\times\catit{B}$-morphisms $(X_1,B_1)\arw{(f,\varphi)}(X_2,B_2)$, making the diagram
 $$
   \xymatrix{B_1 \ar[d]_{\kappa_1} \ar[rr]^{\varphi} & & B_2
             \ar[d]^{\kappa_2}\\
             TX_1 \ar[rr]_{Tf} & & TX_2}
 $$
 commute.
 \sqed
\end{defn}

\par
In the following, we provide some examples of variety-based bornological systems.

\begin{exmp}
 \label{exmp:4}
 \hfill\par
 \begin{enumerate}[(1)]
  \item \bornsyst{(\mcal{P},\latcb)} provides an analogue of the category \cat{BornSys} of bornological systems and bounded morphisms
        of~\cite{Paseka2015}.
  \item \bornsyst{(\mcal{P}_{L}^{\cat{CLat}},\latcb)} provides an analogue of the category \lbornsys{L} of~\cite{Paseka2015}.
  \item \bornsyst{(\mcal{P}_{\cat{S}}^{\cat{CLat}},\latcb)} provides an analogue of the category \lbornsys{\cat{L}}
        of~\cite{Paseka2015}.
        \sqed
 \end{enumerate}
\end{exmp}

\par
We notice that our variety-based approach does not fully incorporate the lattice-valued system setting of~\cite{Paseka2015}. More precisely, in~\cite{Paseka2015}, the maps $\varphi$ and $\kappa_i$, employed in the diagram of Definition~\ref{defn:13}, preserve different algebraic structure, i.e., do not come from the same category. In particular, the definition of $\kappa_i$ is more restrictive than that of $\varphi$ (even in view of Definition~\ref{defn:13}, which puts an additional condition on $\kappa_i$). The next subsection, however, shows that variety-based setting is more convenient for dealing with the so-called spatialization procedure for bornological systems.


\subsection{Bornological systems versus bornological spaces}


\par
In~\cite{Vickers1989}, S.~Vickers showed that the category of topological spaces is isomorphic to a full coreflective subcategory of the category of topological systems. The respective functor $\cat{TopSys}\arw{Spat}\cat{Top}$ (from topological systems to topological spaces) was called the system spatialization procedure. It is the main purpose of this subsection, to show a bornological analogue of the above functor $Spat$.

\begin{prop}
 \label{prop:3}
 There exists a full embedding \incl{\bornt{T}}{E}{\bornsyst{T},} $E((X_1,\tau_1)\arw{f}(X_2,\tau_2))=(X_1,e_{\tau_1},\tau_1)\arw{(f,\overline{Tf})}(X_2,e_{\tau_2},\tau_2)$, where \incl{\tau_i}{e_{\tau_i}}{TX_i} is the embedding, and $\overline{Tf}$ is the restriction $(Tf)|_{\tau_1}^{\tau_2}$.
\end{prop}
\begin{pf}
 The functor is correct on both objects and morphism, which follows from the commutative diagram
 $$
   \xymatrix@1@H=12pt{\tau_1 \ar@{^{(}->}[d]_{e_{\tau_1}} \ar[rr]^{\overline{Tf}=Tf|_{\tau_1}^{\tau_2}} & & \tau_2 \ar@{^{(}->}[d]^{e_{\tau_2}}\\
               TX_1 \ar[rr]_{Tf} & & TX_2.}
 $$
 The same diagram takes care of fullness (every \catit{B}-morphism $\tau_1\arw{\varphi}\tau_2$, which makes the above diagram commute, will be then necessarily the restriction of $Tf$).
 \qed
\end{pf}

\begin{prop}
 \label{prop:4}
 There exists a functor $\bornsyst{T}\arw{Spat}\bornt{T}$ defined by $Spat((X_1,\kappa_1,B_1)\arw{(f,\varphi)}(X_2,\kappa_2,B_2))=(X_1,\tau_1=\gnrt{\mcal{P}\kappa_1(B_1)}_{\cat{B}})\arw{f}
 (X_2,\tau_2=\gnrt{\mcal{P}\kappa_2(B_2)}_{\cat{B}})$.
\end{prop}
\begin{pf}
 To show that $Spat$ is correct on objects, we notice that $\gnrt{\tau_i}_{\cat{A}}=\gnrt{\gnrt{\mcal{P}\kappa(B_i)}_{\cat{B}}}_{\cat{A}}=
 \gnrt{\mcal{P}\kappa(B_i)}_{\cat{A}}\ni\top$.

 \par
 To show that the functor is correct on morphisms, we consider the commutative diagram of Definition~\ref{defn:13}. Given an ideal term $p(\vec{x},\vec{y})$, for every $\vec{x}\in(TX_1)^n$ and every $\vec{b}\in B_1^m$,
 $Tf(p(\vec{x},\kappa_1(\vec{b})))=\prm{p}(\vec{z},Tf\circ\kappa_1(\vec{b}))=
 \prm{p}(\vec{z},\kappa_2\circ\varphi(\vec{b}))\in
 \gnrt{\mcal{P}\kappa_2(B_2)}_{\cat{B}}$, which (in view Proposition~\ref{prop:1}) proves the claim.
 \qed
\end{pf}

\par
The main result of this subsection is then as follows.

\begin{thm}
 \label{thm:3}
 $Spat$ is a left-adjoin-left-inverse to $E$.
\end{thm}
\begin{pf}
 Given a \bornsyst{T}-object $(X,\kappa,B)$, $ESpat(X,\kappa,B)=E(X,\tau=\gnrt{\mcal{P}\kappa(B)}_{\cat{B}})=(X,e_{\tau},\tau)$, and, moreover, the diagram
 $$
   \xymatrix@1@H=12pt{B \ar[d]_-{\kappa} \ar[rr]^-{\overline{\kappa}=\kappa|^{\tau}} & & \tau=\gnrt{\mcal{P}\kappa(B)}_{\cat{B}}
                      \ar@{^{(}->}[d]^-{e_{\tau}}\\
                      TX \ar[rr]_-{1_{TX}} & & TX}
 $$
 commutes. Thus, $(X,\kappa,B)\arw{(1_X,\overline{\kappa})}ESpat(X,\kappa,B)$ is a \bornsyst{T}-morphism.

 \par
 To show that $((1_X,\overline{\kappa}),Spat(X,\kappa,B))$ is an $E$-universal arrow for $(X,\kappa,B)$, we take a \bornsyst{T}-morphism $(X,\kappa,B)\arw{(f,\varphi)}(E(\prm{X},\prm{\tau})=(\prm{X},e_{\prm{t}},\prm{\tau}))$ and show that $(Spat(X,\kappa,B)=(X,\tau=\gnrt{\mcal{P}\kappa(B)}_{\cat{B}}))\arw{f}(\prm{X},\prm{\tau})$ is a \bornt{T}-morphism. Commutativity of the diagram
 $$
   \xymatrix@1@H=12pt{B \ar[d]_-{\kappa} \ar[rr]^-{\varphi} & & \prm{\tau}
                      \ar@{^{(}->}[d]^-{e_{\prm{\tau}}}\\
                      TX \ar[rr]_-{Tf} & & T\prm{X}}
 $$
 gives then that for every ideal term $p(\vec{x},\vec{y})$, every $\vec{x}\in(TX)^n$ and every $\vec{b}\in B^m$, $Tf(p(\vec{x},\kappa(\vec{b})))=\prm{p}(\vec{z},Tf\circ\kappa(\vec{b}))=\prm{p}(\vec{z},e_{\prm{\tau}}\circ\varphi(\vec{b}))=
 \prm{p}(\vec{z},\varphi(\vec{b}))\in\prm{\tau}$, which was to show.

 \par
 It is easy to see that the triangle
 $$
   \xymatrix@1{(X,\kappa,B)\ar[rr]^-{(1_X,\overline{\kappa})} \ar"2,3"_-{(f,\varphi)} & &
               ESpat(X,\kappa,B) \ar@{.>}[d]^-{Ef}\\
               & & E(\prm{X},\prm{\tau})}
 $$
 commutes, and that $f$ is the unique morphism with such property.

 \par
 To show that $Spat$ is a left inverse to $E$, we notice that given a \bornt{T}-object $(X,\tau)$, $SpatE(X,\tau)=Spat(E,e_{\tau},\tau)=(X,\gnrt{\mcal{P}e_{\tau}(\tau)}_{\cat{B}})=(X,\tau)$.
 \qed
\end{pf}

\par
As an immediate consequence of Propositions~\ref{prop:3}, \ref{prop:4} and Theorem~\ref{thm:3}, we get the following result.

\begin{cor}
 \label{cor:1}
 \bornt{T} is isomorphic to a full reflective subcategory of \bornsyst{T}.
\end{cor}

\par
We are currently unable to verify, whether the reflective subcategory of Corollary~\ref{cor:1} is actually epireflective or (regular epi)-reflective, which depends on whether the \catit{B}-morphism $B\arw{\overline{\kappa}}\gnrt{\mcal{P}\kappa(B)}_{\cat{B}}$ of the proof of Theorem~\ref{thm:3} is a (regular) epimorphism. More precisely, we presently lack a proper characterization of (regular) epimorphisms in the category \catit{B}, which will be the subject of our further study.


\section{Conclusion and future work}


\par
In this paper, we introduced a new approach to (lattice-valued) bornology, which is based in varieties of algebras, and which is motivated by variety-based topology of~\cite{Solovyov2012,Solovyov2013a}. Our main idea here (similar to the case of topology) is to provide a common setting for different lattice-valued approaches to bornology. In particular, our variety-based approach incorporates the settings of M.~Abel and A.~\v{S}ostak~\cite{Abel2011} as well as our previous lattice-valued approach of~\cite{Paseka2015}. Moreover, in this new setting, it is possible to get a bornological analogue of the notion of topological system of S.~Vickers~\cite{Vickers1989}, and show that the category of bornological spaces is isomorphic to a full reflective subcategory of the category of bornological systems. We would like to end the paper with several open problems, which concern its topic of study.


\subsection{Topological nature of the category of bornological spaces}


\par
In Theorem~\ref{thm:1} of this paper, we showed sufficient conditions for the category \bornt{T} of variety-based bornological spaces to be topological over its ground category \cat{X}. Moreover, in~\cite{Paseka2015}, we showed the necessary and sufficient conditions for the category \lborn{\cat{L}} of lattice-valued variable-basis bornological spaces to be topological. The first open problem of this paper can be formulated then as follows.

\begin{prob}
 \label{prob:1}
 What are the necessary and sufficient conditions for the category \bornt{T} to be topological over its ground category \cat{X}?
 \sqed
\end{prob}


\subsection{The nature of the category of bornological systems}


\par
In~\cite{Solovyov2012}, we showed that the category of variety-based topological systems is essentially algebraic (in the sense of~\cite[Definition~23.5]{Adamek2009}), which, taken together with the fact that the topological category of variety-based topological spaces is isomorphic to a full coreflective subcategory of the category of variety-based topological systems, provides then an embedding of topology into algebra~\cite{Solovjovs2009}. The next open problem of this paper can be formulated therefore as follows.

\begin{prob}
 \label{prob:2}
 What is the nature of the category \bornsyst{T} of variety-based bornological systems?
 \sqed
\end{prob}


\subsection{Localification procedure for bornological systems}


\par
In~\cite{Vickers1989}, S.~Vickers provided two procedures for topological systems, namely, spatialization and localification. The former has already been considered in this paper, whereas the latter not. Briefly speaking, the category of locales is isomorphic to a full reflective subcategory of the category of topological systems. The functor $\cat{TopSys}\arw{Loc}\cat{Loc}$ (from topological systems to locales) is then called localification procedure for topological systems. In \cite{Solovyovi}, we provided a variety-based generalization of the localification procedure. The last open problem of this paper is then as follows.

\begin{prob}
 \label{prob:3}
 Provide a variety-based localification procedure for bornological systems.
 \sqed
\end{prob}

\par
We notice that one part of the above procedure is almost trivial.

\begin{prop}
 \label{prop:5}
 There exists a functor $\bornsyst{T}\arw{Loc}\catit{B}$, which is defined by $Loc((X_1,\kappa_1,B_1)\arw{(f,\varphi)}(X_2,\kappa_2,B_2))=
 B_1\arw{\varphi}B_2$.
\end{prop}

\par
It is the functor in the opposite direction, which causes the main problem. We also notice that the term ``localification" itself, which stems from ``locale" in the topological setting, is not quite correct in case or bornology. In~\cite{Paseka2015}, we have proposed the term \emph{bornale} for the underlying algebraic structures of lattice-valued bornology. Thus, the respective procedure could be called then ``bornalification".

\par
The above open problems will be addressed in our next papers on the topic of variety-based bornology.


\section*{Acknowledgements}


\par
This is a pre-print of an article published in Fuzzy Sets and Systems. 
The final authenticated version of the article is available online at: \newline 
https://www.sciencedirect.com/science/article/pii/S0165011415003267.

\end{document}